\documentclass[12pt,a4paper,reqno]{amsart}
\usepackage{fullpage}
\usepackage{verbatim}
\usepackage{hyperref}
\usepackage{pslatex}

\newcommand{\odip}[2]{o _{#1}\!\left(#2\right)\mathchoice{\!}{}{}{}}
\newcommand{\odi}[1]{\odip{}{#1}}
\newcommand{\Odip}[2]{\mathcal{O}_{#1}\left(#2\right)}
\newcommand{\Odi}[1]{\mathcal{O}\left(#1\right)}

\newcommand{\dx}{\mathrm{d}}
\newcommand{\e}{\mathrm{e}}
\newcommand{\ii}{\mathrm{i}}
\newcommand{\E}{\widetilde{\mathcal{E}}}
\newcommand{\C}{\mathcal{C}}
\newcommand{\eps}{\varepsilon}
\newcommand{\Stilde}{\widetilde{S}}
\newcommand{\K}{\mathbf{k}}
\newcommand{\gA}{\mathfrak{A}}
\newcommand{\gB}{\mathfrak{B}}

\newtheorem{Theorem}{Theorem}[section]
\newtheorem{Lemma}{Lemma}[section]

\allowdisplaybreaks

\title{A note on an average additive problem with prime numbers}

\author[M.~Cantarini, A.~Gambini, A.~Zaccagnini]
       {Marco Cantarini, Alessandro Gambini, Alessandro Zaccagnini}

\date{\today}


\begin{document}

\begin{abstract}
We continue our investigations on the average number of representations
of a large positive integer as a sum of given powers of prime
numbers.
The average is taken over a ``short'' interval, whose admissible length
depends on whether or not we assume the Riemann Hypothesis.
\end{abstract}

\maketitle

\section{Introduction}

We pursue recent investigations by the present authors and
Alessandro Languasco in \cite{CantariniGZ2018} and
\cite{CantariniGLZ2018}.
In this short note we study general average additive problems: let
$\K = (k_1, k_2, \dots, k_r)$, where $2 \le k_1 \le k_2 \le \dots \le k_r$
and $k_j$ is an integer for all $j \in \{1, \dots, r \}$.
Let
\begin{equation}
\label{Rk-def}
  R(n; \K)
  =
  \sum_{n = m_1^{k_1} + \cdots + m_r^{k_r}}
    \Lambda(m_1) \cdots \Lambda(m_r),
\end{equation}
where $\Lambda$ is the von Mangoldt function, that is,
$\Lambda(p^m) = \log(p)$ if $p$ is a prime number and $m$ is a
positive integer, and $\Lambda(n) = 0$ for all other integers.
We write $\rho = \rho(\K) = k_1^{-1} + \cdots + k_r^{-1}$, for the
``density'' of the problem, $\gamma_k = \Gamma(1 + 1 / k)$ where
$\Gamma$ is the Euler Gamma-function and
$G(\K) = \gamma_{k_1} \cdots \gamma_{k_r}$.

Proving the expected individual asymptotic formula for $R(n; \K)$ as
$n \to \infty$ along ``admissible'' residue classes (that is, avoiding
those residue classes which can not contain values of the form
$p_1^{k_1} + \cdots + p_r^{k_r}$ because of the uneven distribution
of prime powers in residue classes) is very difficult if either $r$ or
$\rho$ is small.
Our main goal is to give an asymptotic formula for the average value
of $R(n; \K)$ for $n \in [N + 1, N + H]$ where $N \to +\infty$
and $H = H(N; \K)$ is as small as possible.
Here we assume $r \ge 3$, since binary problems have been thoroughly
studied in \cite{LanguascoZ2016c}, \cite{LanguascoZ2017b},
\cite{LanguascoZ2018b}, \cite{LanguascoZ2018a}.

\begin{Theorem}
\label{Th-unconditional}
Let $\K = (k_1, \dots, k_r)$, where $2 \le k_1 \le \dots \le k_r$, be an
$r$-tuple of integers with $r \ge 3$.
For every $\eps > 0$ there exists a constant $C = C(\eps) > 0$,
independent of $\K$, such that
\[
  \sum_{n = N + 1}^{N + H} R(n; \K)
  =
  \frac{G(\K)}{\Gamma(\rho)}
  H N^{\rho - 1}
  +
  \Odip{\K}{H N^{\rho - 1}
    \exp\Bigl\{ -C \Bigl(\frac{\log N}{\log\log N}\Bigr)^{1/3} \Bigr\}}
\]
as $N \to +\infty$, uniformly for
$N^{1 - 5 / (6 k_r) + \eps} < H < N^{1 - \eps}$.
\end{Theorem}

It is well known that the Riemann Hypothesis (RH for short) implies
that prime numbers are fairly regularly distributed.
In this problem, it has the effect of allowing far wider ranges for
$H$, that is, much smaller values of $H$ are admissible than in
Theorem~\ref{Th-unconditional}.
The final error term is also smaller, as it is to be expected.

We use throughout the paper the convenient notation $f = \infty(g)$
as equivalent to $g = \odi{f}$.

\begin{Theorem}
\label{Th-RH}
Assume the Riemann Hypothesis.
Let $\K = (k_1, \dots, k_r)$, where $2 \le k_1 \le \dots \le k_r$, be an
$r$-tuple of integers with $r \ge 3$.
For every $\eps > 0$ there exists a constant $C = C(\eps) > 0$,
independent of $\K$, such that
\[
  \sum_{n = N + 1}^{N + H} R(n; \K)
  =
  \frac{G(\K)}{\Gamma(\rho)}
  H N^{\rho - 1}
  +
  \Odip{\K}{H^2 N^{\rho - 2} + H^{1 / 2} N^{\rho - 1/2 - 1 / (2 k_r)} L^3}
\]
as $N \to +\infty$, uniformly for
$H = \infty\bigl( N^{1 - 1 / k_r} (\log N)^6 \bigr)$ with
$H < N^{1 - \eps}$.
\end{Theorem}

Theorem~\ref{Th-unconditional} contains as special cases all results
in \cite{CantariniGZ2018} and \cite{CantariniGLZ2018}, whereas
Theorem~\ref{Th-RH} is occasionally slightly weaker because our basic
combinatorial identity here, equation~\eqref{comb-identity}, is less
efficient than the identities we used in the papers mentioned above.

\section{Definitions and preparation for the proofs}

We rewrite $R(N; \K)$ as the integral over the unit interval of
the product of suitable exponential sums.
We then proceed to ``replace'' each exponential sum by its
approximation, which is given by the leading term of the Prime Number
Theorem.
This gives rise to the main term and also to a number of additional
terms that we have to bound in various ways.
Let $S_j = x_j + y_j$ for $j \in \{1$, \dots, $r \}$: then we have
\begin{equation}
\label{comb-identity}
  \prod_{j = 1}^r S_j
  =
  \prod_{j = 1}^r (x_j + y_j)
  =
  \prod_{j = 1}^r x_j + \gA + \gB,
\end{equation}
where
\begin{align}
\label{def-gA}
  \gA
  &=
  \sum_{i = 1}^r y_i \Bigl( \prod_{j \ne i} S_j \Bigr), \\
\label{def-gB}
  \gB
  &=
  \sum_{\substack{I \subseteq \{ 1, \dots, r \} \\ |I| \ge 2}}
    c_r(I)
    \Bigl( \prod_{i \in \{ 1, \dots, r \} \setminus I} x_i \Bigr)
    \Bigl( \prod_{i \in I} y_i \Bigr),
\end{align}
for suitable coefficients $c_r(I)$.
In fact, according to the definitions~\eqref{def-Stilde}
and~\eqref{def-Etilde} below, we will choose
$\Stilde_{k_j}(\alpha) = S_j = x_j + y_j$ where
$x_j = x_j(\alpha) = \gamma_j z^{-1/k_j}$ and
$y_j = y_j(\alpha) = \E_{k_j}(\alpha)$, so that we can exploit the
fact that $S_j, x_j, y_j \ll N^{1 / k_j}$ and that $y_j$ is small in
$L^2$-norm by Lemma~\ref{LP-Lemma-gen} below.

For real $\alpha$ we write $\e(\alpha) = \e^{2 \pi \ii \alpha}$.
We take $N$ as a large positive integer, and write $L = \log N$ for
brevity.
In this and in the following section $k$ denotes any positive real
number.
Let $z = 1 / N - 2 \pi \ii \alpha$ and
\begin{equation}
\label{def-Stilde}
  \Stilde_k(\alpha)
  =
  \sum_{n \ge 1} \Lambda(n) \e^{-n^k / N} \e(n^k \alpha)
  =
  \sum_{n \ge 1} \Lambda(n) \e^{- n^k z}.
\end{equation}
Thus, recalling definition \eqref{Rk-def} and using \eqref{def-Stilde},
for all $n \ge 1$ we have
\begin{equation}
\label{basic-Rk}
  R(n; \K)
  =
  \sum_{n_1^{k_1} + \cdots + n_r^{k_r} = n}
    \Lambda(n_1) \cdots \Lambda(n_r)
  =
  \e^{n / N}
  \int_{-1/2}^{1/2}
    \Stilde_{k_1}(\alpha) \cdots \Stilde_{k_r}(\alpha)
    \, \e(-n \alpha) \, \dx \alpha.
\end{equation}
It is clear from the above identity that we are only interested in the
range $\alpha \in [-1/2, 1/2]$.
We record here the basic inequality
\begin{equation}
\label{z-bound}
  \vert z \vert^{-1}
  \ll
  \min \{ N, \vert \alpha \vert^{-1} \}.
\end{equation}
We also need the following exponential sum over the ``short interval''
$[1, H]$
\[
  U(\alpha, H)
  =
  \sum_{m = 1}^H \e(m \alpha),
\]
where $H \le N$ is a large integer.
We recall the simple inequality
\begin{equation}
\label{U-bound}
  \vert U(\alpha, H) \vert
  \le
  \min \{ H, \vert \alpha \vert^{-1} \}.
\end{equation}
With these definitions in mind and recalling \eqref{basic-Rk}, our
starting point is the identity
\begin{equation}
\label{basic-identity}
  \sum_{n = N + 1}^{N + H}
    \e^{-n / N} R(n; \K)
  =
  \int_{-1/2}^{1/2}
    \Stilde_{k_1}(\alpha) \cdots \Stilde_{k_r}(\alpha)
    U(-\alpha, H) \, \e(-N \alpha) \, \dx \alpha.
\end{equation}
The basic strategy is to replace each factor $\Stilde_k(\alpha)$ by
its expected main term, which is $\gamma_k / z^{1 / k}$, and
estimating the ensuing error term by means of a combination of
techniques and bounds for exponential sums, with the aid
of~\eqref{comb-identity}.
One key ingredient is the $L^2$-bound in
Lemma~\ref{LP-Lemma-gen}, which we may use only in a restricted
range, and we need a different argument on the remaining part of the
integration interval.
This leads to some complications in details.
The conditional case, when the Riemann Hypothesis is assumed, has
a somewhat simpler proof, as we see in \S\ref{proof-cond}.

\section{Lemmas}

For brevity, we set
\begin{equation}
\label{def-Etilde}
  \E_k(\alpha)
  :=
  \Stilde_k(\alpha)
  -
  \frac{\gamma_k}{z^{1 / k}}
  \qquad\text{and}\qquad
  A(N; c)
  :=
  \exp\Bigl\{ c \Bigl(\frac{\log N}{\log\log N}\Bigr)^{1/3} \Bigr\},
\end{equation}
where $c$ is a real constant.

\begin{Lemma}[Lemma~3 of \cite{LanguascoZ2016a}]
\label{LP-Lemma-gen}
Let $\eps$ be an arbitrarily small positive constant, $k \ge 1$ be
an integer, $N$ be a sufficiently large integer and $L = \log N$.
Then there exists a positive constant $c_1 = c_1(\eps)$, which does
not depend on $k$, such that
\[
  \int_{-\xi}^{\xi}
    \bigl\vert \E_k(\alpha) \bigr\vert^2 \, \dx \alpha
  \ll_k
  N^{2 / k - 1} A(N; - c_1)
\]
uniformly for $0 \le \xi < N^{ -1 + 5 / (6 k) - \eps}$.
Assuming the Riemann Hypothesis we have
\[
  \int_{-\xi}^{\xi} \,
    \bigl\vert \E_k(\alpha) \bigr\vert^2 \, \dx \alpha
  \ll_k
  N^{1 / k}\xi L^2
\]
uniformly for $0 \le \xi \le 1 / 2$.
\end{Lemma}

We remark that the proof of Lemma~3 in \cite{LanguascoZ2016a} contains
oversights which are corrected in \cite{LanguascoZ2017e}.
The next result is a variant of Lemma~4 of \cite{LanguascoZ2016a},
which is fully proved in \cite{CantariniGZ2018}.

\begin{Lemma}
\label{Laplace-formula}
Let $N$ be a positive integer, $z = z(\alpha) = 1 / N - 2 \pi \ii \alpha$,
and $\mu > 0$.
Then, uniformly for $n \ge 1$ and $X > 0$ we have
\[
  \int_{-X}^X z^{-\mu} \e(-n \alpha) \, \dx \alpha
  =
  \e^{- n / N} \frac{n^{\mu - 1}}{\Gamma(\mu)}
  +
  \Odip{\mu}{\frac1{n X^{\mu}}}.
\]
\end{Lemma}

\begin{Lemma}[Lemma 3.3 of \cite{CantariniGZ2018}]
\label{Stilde-bound}
We have $\Stilde_k(\alpha) \ll_k N^{1 / k}$.
\end{Lemma}

We record an immediate consequence of \eqref{z-bound},
\eqref{def-Etilde} and Lemma~\ref{Stilde-bound}:
\begin{equation}
\label{E-bound}
  \E_k(\alpha)
  \ll_k
  N^{1 / k}.
\end{equation}
Our next tool is the extension to $\Stilde_k$ of Lemma~7 of Tolev
\cite{Tolev1992}.
The proof can be found in \cite{CantariniGLZ2018}.

\begin{Lemma}
\label{Tolev-Lemma}
Let $k > 1$ and $\tau > 0$. Then
\[
  \int_{-\tau}^{\tau} \vert \Stilde_k(\alpha) \vert^2 \, \dx \alpha
  \ll
  \bigl(\tau N^{1/k} + N^{2/k - 1}\bigr) L^3.
\]
\end{Lemma}

\begin{Lemma}[Lemma 3.6 of \cite{CantariniGZ2018}]
\label{mt-evaluation}
For $N \to +\infty$, $H \in [1, N]$ and a real number $\lambda$ we
have
\[
  \sum_{n = N + 1}^{N + H}
    \e^{- n / N} n^{\lambda}
  =
  \frac1{\e}
  H N^{\lambda}
  +
  \Odip{\lambda}{H^2 N^{\lambda - 1}}.
\]
\end{Lemma}

\section{Proof of Theorem \ref{Th-unconditional}}

We need to introduce another parameter $B = B(N)$, defined as
\begin{equation}
\label{def-B}
  B = N^{2 \eps}.
\end{equation}
We can not take $B = 1$, because of the estimate in \S\ref{sub-I4}.
We let $\C = \C(B,H) = [-1/2, -B/H] \cup [B/H, 1/2]$, and
write $\Stilde_{k_j}(\alpha) = x_j + y_j$ where
$x_j = x_j(\alpha) = \gamma_j z^{-1/k_j}$ and
$y_j = y_j(\alpha) = \E_{k_j}(\alpha)$ in~\eqref{comb-identity}, so
that
\begin{equation}
\label{identity}
  \Stilde_{k_1}(\alpha) \cdots \Stilde_{k_r}(\alpha)
  =
  \prod_{j = 1}^r (x_j + y_j)
  =
  \prod_{j = 1}^r x_j(\alpha)
  +
  \gA(\alpha) + \gB(\alpha),
\end{equation}
where $\gA(\alpha)$ and $\gB(\alpha)$ are defined by~\eqref{def-gA}
and~\eqref{def-gB} respectively.
We multiply \eqref{identity} by $U(-\alpha, H) \* \e(-N \alpha)$
and integrate over the interval $[-B / H, B / H]$.
Recalling \eqref{basic-identity} we have
\begin{align*}
  \sum_{n = N + 1}^{N + H}
    \e^{-n / N} R(n; \K)
  &=
  G(\K)
  \int_{-B/H}^{B/H} \frac{U(-\alpha, H)}{z^{\rho}}
    \, \e(-N \alpha) \, \dx \alpha \\
  &\qquad+
  \int_{-B/H}^{B/H} \gA(\alpha) U(-\alpha, H) \e(-N \alpha) \, \dx \alpha \\
  &\qquad+
  \int_{-B/H}^{B/H} \gB(\alpha) U(-\alpha, H) \e(-N \alpha) \, \dx \alpha \\
  &\qquad+
  \int_\C \Stilde_{k_1}(\alpha) \cdots \Stilde_{k_r}(\alpha)
    U(-\alpha, H) \e(-N \alpha) \, \dx \alpha \\
  &=
  G(\K) I_1
  +
  I_2 + I_3 + I_4,
\end{align*}
say.
The first summand gives rise to the main term via
Lemma~\ref{Laplace-formula}, the next two are majorised in
\S\ref{sub-I2}--\ref{sub-I3} by means of Lemma~\ref{Stilde-bound} and
the $L^2$-estimate provided by Lemma~\ref{LP-Lemma-gen}.
Finally, $I_4$ is easy to bound using Lemma~\ref{Tolev-Lemma}.

\subsection{Evaluation of \texorpdfstring{$I_1$}{I1}}

It is a straightforward application of Lemma~\ref{Laplace-formula}:
we have
\begin{equation}
\label{prep-mt}
  \int_{-B/H}^{B/H} \frac{U(-\alpha, H)}{z^{\rho}}
    \, \e(-N \alpha) \, \dx \alpha
  =
  \frac1{\Gamma(\rho)}
  \sum_{n = N + 1}^{N + H}
    \e^{- n / N} n^{\rho - 1}
  +
  \Odip{\K}{\frac HN \Bigl( \frac HB \Bigr)^{\rho}}.
\end{equation}
We evaluate the sum on the right-hand side of \eqref{prep-mt} by means
of Lemma~\ref{mt-evaluation} with $\lambda = \rho - 1$.
Summing up, we have
\begin{equation}
\label{final-mt}
  \int_{-B/H}^{B/H} \frac{U(-\alpha, H)}{z^{\rho}}
    \, \e(-N \alpha) \, \dx \alpha
  =
  \frac1{\e \Gamma(\rho)} H N^{\rho - 1}
  +
  \Odip{\K}{H^2 N^{\rho - 2} + \frac HN \Bigl( \frac HB \Bigr)^{\rho}}.
\end{equation}

We now choose the range for $H$: since will need
Lemma~\ref{LP-Lemma-gen}, we see that we can take
\begin{equation}
\label{H-bound}
  H > N^{1 - 5 / (6 k_r) + 3 \eps}.
\end{equation}

\subsection{Bound for \texorpdfstring{$I_2$}{I2}}
\label{sub-I2}

We recall the bound~\eqref{U-bound}, and Lemmas~\ref{Stilde-bound} and
\ref{Tolev-Lemma}.
Using Lemma~\ref{LP-Lemma-gen} and the Cauchy-Schwarz inequality where
appropriate, we see that the contribution from
$\Stilde_{k_1}(\alpha) \* \cdots \* \Stilde_{k_{r-1}}(\alpha) y_r$,
say, is
\begin{align}
\notag
  &\ll_{\K}
  H
  \max_{\alpha \in [-1/2, 1/2]}
    \vert \Stilde_{k_1}(\alpha) \cdots \Stilde_{k_{r-2}(\alpha)} \vert
  \Bigl(
  \int_{-B / H}^{B / H} \vert \Stilde_{k_{r-1}}(\alpha) \vert^2 \, \dx \alpha
  \int_{-B / H}^{B / H} \vert \E_{k_r}         (\alpha) \vert^2 \, \dx \alpha
  \Bigr)^{1 / 2} \\
\notag
  &\ll_{\K}
  H N^{1 / k_1 + \cdots + 1 / k_{r-2}} L^{3 / 2}
  \Bigl( \frac BH N^{1 / k_{r-1}} + N^{2 / k_{r-1} - 1} \Bigr)^{1 / 2}
  \bigl( N^{2 / k_r - 1} A(N; - c_1) \bigr)^{1 / 2} \\
\label{bound-I2}
  &\ll_{\K}
  H N^{\rho - 1} A \Bigl(N; - \frac13 c_1\Bigr),
\end{align}
where $c_1 = c_1(\eps) > 0$ is the constant provided by
Lemma~\ref{LP-Lemma-gen}, which we can use on the interval
$[-B / H, B / H]$ since $B$ and $H$ satisfy \eqref{def-B} and
\eqref{H-bound} respectively.
The other summands in $I_2$ are treated in the same way.

\subsection{Bound for \texorpdfstring{$I_3$}{I3}}
\label{sub-I3}

We remark that, by definition~\eqref{def-gB}, each summand
in $\gB(\alpha)$ is the product of $r$ factors chosen among
the $x_j$s and the $y_j$s, with at least two of the latter type.
Using \eqref{z-bound}, \eqref{U-bound} and Lemma~\ref{LP-Lemma-gen},
by the Cauchy-Schwarz inequality, we see that the contribution
from the term $y_1 y_2 x_3 \dots x_r$, say, is
\begin{align}
\notag
  &=
  \gamma_{k_3} \cdots \gamma_{k_r}
  \int_{-B/H}^{B/H}
    \frac{\E_{k_1}(\alpha) \E_{k_2}(\alpha)}{z^{1/k_3 + \cdots + 1/k_r}}
    U(-\alpha, H) \e(-N \alpha) \, \dx \alpha \\
\notag
  &\ll_{\K}
  H N^{1 / k_3 + \cdots + 1 / k_r}
  \Bigl(
    \int_{-B/H}^{B/H} \vert \E_{k_1}(\alpha) \vert^2 \, \dx \alpha
    \int_{-B/H}^{B/H} \vert \E_{k_2}(\alpha) \vert^2 \, \dx \alpha
  \Bigr)^{1/2} \\
\label{bound-I3}
  &\ll_{\K}
  H N^{\rho - 1} A(N; -c_1).
\end{align}
Furthermore, we recall the bound $\E_k(\alpha) \ll_k N^{1 / k}$
in \eqref{E-bound}.
Hence we may treat the other summands in $I_3$ in the same way, since
$x_j$, $y_j \ll_{k_j} N^{1 / k_j}$ for $j \in \{ 1, \dots, r \}$.

\subsection{Bound for \texorpdfstring{$I_4$}{I4}}
\label{sub-I4}

Using a partial integration from Lemma~\ref{Tolev-Lemma} and the
Cauchy-Schwarz inequality, we have
\begin{align}
\notag
  I_4
  &=
  \int_\C \Stilde_{k_1}(\alpha) \cdots \Stilde_{k_r}(\alpha)
    U(-\alpha, H) \e(-N \alpha) \, \dx \alpha \\
\notag
  &\ll_{\K}
  \max_{\alpha \in [-1/2, 1/2]}
    \vert \Stilde_{k_1}(\alpha) \cdots \Stilde_{k_{r - 2}}(\alpha) \vert
  \Bigl(
    \int_\C \vert \Stilde_{k_{r-1}}(\alpha) \vert^2
      \, \frac{\dx \alpha}{\vert \alpha \vert}
    \int_\C \vert \Stilde_{k_r}    (\alpha) \vert^2
      \, \frac{\dx \alpha}{\vert \alpha \vert}
  \Bigr)^{1/2} \\
\label{bound-I4}
  &\ll_{\K}
  N^{1 / k_1 + \dots + 1 / k_{r-2}}
  \Bigl( \frac{H^2}{B^2} N^{2 / k_{r-1} + 2 / k_r - 2} L^6 \Bigr)^{1/2}
  \ll_{\K}
  \frac HB N^{\rho - 1} L^3,
\end{align}
because of \eqref{H-bound}.
This is $\ll_{\K} H N^{\rho - 1} A(N; -c_1 / 3)$, by our choice
in~\eqref{def-B}.

\subsection{Completion of the proof}
\label{final-Th1}

For simplicity, from now on we assume that $H \le N^{1 - \eps}$.
Summing up from \eqref{final-mt}, \eqref{bound-I2}, \eqref{bound-I3}
and \eqref{bound-I4}, we proved that
\begin{equation}
\label{smooth-Th1}
  \sum_{n = N + 1}^{N + H}
    \e^{-n / N} R(n; \K)
  =
  \frac{G(\K)}{\e \Gamma(\rho)}
  H N^{\rho - 1}
  +
  \Odip{\K}{H N^{\rho - 1} A \Bigl(N; - \frac13 c_1\Bigr)},
\end{equation}
provided that \eqref{def-B} and \eqref{H-bound} hold, since the other
error terms are smaller in our range for $H$.
In order to achieve the proof, we have to remove the exponential
factor on the left-hand side, exploiting the fact that, since $H$ is
``small,'' it does not vary too much over the summation range.
Since $\e^{-n / N} \in [\e^{-2}, \e^{-1}]$ for all
$n \in [N + 1, N + H]$,
we can easily deduce from \eqref{smooth-Th1} that
\[
  \e^{-2}
  \sum_{n = N + 1}^{N + H} R(n; \K)
  \le
  \sum_{n = N + 1}^{N + H}
    \e^{-n / N} R(n; \K)
  \ll_{\K}
  H N^{\rho - 1}.
\]
We can use this weak upper bound to majorise the error term arising
from the development $\e^{-x} = 1 + \Odi{x}$ that we need in the
left-hand side of \eqref{smooth-Th1}.
In fact, we have
\begin{align*}
  \sum_{n = N + 1}^{N + H}
    \e^{-n / N} R(n; \K)
  &=
  \sum_{n = N + 1}^{N + H}
    \bigl(\e^{-1} + \Odi{(n - N) N^{-1}} \bigr) R(n; \K) \\
  &=
  \e^{-1}
  \sum_{n = N + 1}^{N + H} R(n; \K)
  +
  \Odip{\K}{ H^2 N^{\rho - 2} }.
\end{align*}
Finally, substituting back into \eqref{smooth-Th1}, we obtain the
required asymptotic formula for $H$ as in the statement of
Theorem~\ref{Th-unconditional}.

\section{Proof of Theorem \ref{Th-RH}}
\label{proof-cond}

Here we assume the Riemann Hypothesis: as we mentioned above, we
obtain stronger results (wider ranges for $H$, better error term) and
the proof is simpler because Lemma~\ref{LP-Lemma-gen} applies to the
whole unit interval.
In fact, we use identity~\eqref{identity} over $[-1/2, 1/2]$.
Recalling \eqref{basic-identity} we have
\begin{align*}
  \sum_{n = N + 1}^{N + H}
    \e^{-n / N} R(n; \K)
  &=
  G(\K)
  \int_{-1 / 2}^{1 / 2} \frac{U(-\alpha, H)}{z^{\rho}}
    \, \e(-N \alpha) \, \dx \alpha \\
  &\qquad+
  \int_{-1 / 2}^{1 / 2} \gA(\alpha) U(-\alpha, H) \e(-N \alpha) \, \dx \alpha \\
  &\qquad+
  \int_{-1 / 2}^{1 / 2} \gB(\alpha) U(-\alpha, H) \e(-N \alpha) \, \dx \alpha \\
  &=
  G(\K) I_1 + I_2 + I_3,
\end{align*}
say.
For the main term we use Lemma~\ref{Laplace-formula} over $[-1/2, 1/2]$
and then Lemma~\ref{mt-evaluation} with $\lambda = \rho - 1$, obtaining
\begin{equation}
\label{final-mt-rh}
  \int_{-1 / 2}^{1 / 2} \frac{U(-\alpha, H)}{z^{\rho}}
    \, \e(-N \alpha) \, \dx \alpha
  =
  \frac1{\e \Gamma(\rho)} H N^{\rho - 1}
  +
  \Odip{\K}{H^2 N^{\rho - 2} + \frac HN}.
\end{equation}
For the other terms, we split the integration range at $1 / H$.
We use Lemma~\ref{LP-Lemma-gen} and \eqref{U-bound} on the interval
$[-1 / H, 1 / H]$,  and a partial-integration argument
from Lemma~\ref{LP-Lemma-gen} in the remaining range.
In view of future constraints (see \eqref{bound-H-RH} below) we assume
that
\begin{equation}
\label{first-bound-H}
  H
  \ge
  N^{1 - 1 / k_r} L.
\end{equation}
We start bounding the contribution of the term
$\Stilde_{k_1}(\alpha) \* \cdots \* \Stilde_{k_{r-1}}(\alpha) y_r$ in
$\gA(\alpha)$ over $[-1 / H, 1 / H]$.
We have that it is
\begin{align}
\notag
  &\ll_{\K}
  H
  \max_{\alpha \in [-1/2, 1/2]}
    \vert \Stilde_{k_1}(\alpha) \cdots \Stilde_{k_{r-2}}(\alpha) \vert
  \Bigl(
  \int_{-1 / H}^{1 / H} \vert \Stilde_{k_{r-1}}(\alpha) \vert^2 \, \dx \alpha
  \int_{-1 / H}^{1 / H} \vert \E_{k_r}         (\alpha) \vert^2 \, \dx \alpha
  \Bigr)^{1 / 2} \\
\notag
  &\ll_{\K}
  H N^{1 / k_1 + \cdots + 1 / k_{r-2}} L^{3 / 2}
  \Bigl( \frac 1H N^{1 / k_{r-1}} + N^{2 / k_{r-1} - 1} \Bigr)^{1 / 2}
  \bigl( N^{1 / k_r} H^{-1} L^2 \bigr)^{1 / 2} \\
\label{bound-A-rh}
  &\ll_{\K}
  H^{1 / 2} N^{\rho - 1/2 - 1 / (2 k_r)} L^{5/2},
\end{align}
by Lemma~\ref{LP-Lemma-gen}, since we assumed \eqref{first-bound-H}.
The same bound holds for other summands in $I_2$.
As above, we remark that $\gB$ is a finite sum of summands which are
products of $x_j$s and $y_j$s, with at least two factors of the latter
type.
For example, we bound the contribution from the term
$x_1 \dots x_{r-2} y_{r-1} y_r$ in $\gB(\alpha)$ on the same interval:
it is
\begin{align}
\notag
  &=
  \gamma_{k_1} \cdots \gamma_{k_{r-2}}
  \int_{-1/H}^{1/H}
    \frac{\E_{k_{r-1}}(\alpha) \E_{k_r}(\alpha)}
         {z^{1/k_1 + \cdots + 1/k_{r-2}}}
    U(-\alpha, H) \e(-N \alpha) \, \dx \alpha \\
\notag
  &\ll_{\K}
  H N^{1 / k_1 + \cdots + 1 / k_{r-2}}
  \Bigl(
    \int_{-1/H}^{1/H} \vert \E_{k_{r-1}}(\alpha) \vert^2 \, \dx \alpha
    \int_{-1/H}^{1/H} \vert \E_{k_r}    (\alpha) \vert^2 \, \dx \alpha
  \Bigr)^{1/2} \\
\label{bound-B-rh}
  &\ll_{\K}
  H N^{1 / k_1 + \cdots + 1 / k_{r-2}}
  \Bigl( N^{1 / k_{r-1} + 1 / k_r} \frac 1{H^2} L^4 \Bigr)^{1/2}
  \ll_{\K}
  N^{\rho - 1 / (2 k_{r-1}) - 1 / (2 k_r)} L^2.
\end{align}
The other summands in $I_3$ can be treated in the same way,
by~\eqref{E-bound}.

We now deal with  the remaining range $[-1/2,1/2] \setminus [-1/H,1/H]$:
by symmetry, it is enough to treat the interval $[1 / H, 1 / 2]$. 
Arguing as in (16) of \cite{CantariniGZ2018} by partial integration
from Lemma~\ref{LP-Lemma-gen}, for $k > 1$ we have
\begin{equation}
\label{part-int-periph}
  \int_{1 / H}^{1 / 2}
    \bigl\vert \E_k(\alpha) \bigr\vert^2 \, \frac{\dx \alpha}{\alpha}
  \ll_k
  N^{1 / k} L^3.
\end{equation}
A partial integration from Lemma~\ref{Tolev-Lemma} also yields
\begin{equation}
\label{Tolev-appl}
  \int_{1 / H}^{1 / 2}
    \vert \Stilde_k(\alpha) \vert^2 \, \frac{\dx \alpha}\alpha
  \ll_k
  N^{1 / k} L^4 + H N^{(2 - k) / k} L^3.
\end{equation}
Proceeding as above, we start bounding the contribution of the term
$\Stilde_{k_1}(\alpha) \* \cdots \* \Stilde_{k_{r-1}}(\alpha) y_r$ in
$\gA(\alpha)$ over $[-1/2,1/2] \setminus [-1 / H, 1 / H]$.
We have that it is
\begin{align}
\notag
  &\ll_{\K}
  \max_{\alpha \in [-1/2, 1/2]}
    \vert \Stilde_{k_1}(\alpha) \cdots \Stilde_{k_{r-2}}(\alpha) \vert
  \Bigl(
  \int_{1 / H}^{1 / 2}
    \vert \Stilde_{k_{r-1}}(\alpha) \vert^2 \, \frac{\dx \alpha}\alpha
  \int_{1 / H}^{1 / 2}
    \vert \E_{k_r}(\alpha) \vert^2 \, \frac{\dx \alpha}\alpha
  \Bigr)^{1 / 2} \\
\notag
  &\ll_{\K}
  N^{1 / k_1 + \cdots + 1 / k_{r-2}}
  \Bigl( N^{1 / k_{r-1}} L^4 + H N^{(2 - k_{r-1}) / k_{r-1}} L^3 \Bigr)^{1 / 2}
  \bigl( N^{1 / k_r}     L^3 \bigr)^{1 / 2} \\
\label{bound-A-rh-minor}
  &\ll_{\K}
  H^{1 / 2} N^{\rho - 1/2 - 1 / (2 k_r)} L^3,
\end{align}
since we assumed \eqref{first-bound-H}.
The other summands in $I_2$ can be estimated in the same way.
Finally, we bound the contribution from the term
$x_1 \dots x_{r-2} y_{r-1} y_r$ in $\gB(\alpha)$ on the same interval:
this is enough in view of our remarks above.
By \eqref{part-int-periph} we may say that it is
\begin{align}
\notag
  &\ll_{\K}
  \int_{1/H}^{1/2}
    \frac{\vert \E_{k_{r-1}}(\alpha) \E_{k_r}(\alpha) \vert}
         {\vert z \vert^{1/k_1 + \cdots + 1/k_{r-2}}}
    \, \frac{\dx \alpha}\alpha \\
\notag
  &\ll_{\K}
  N^{1 / k_1 + \cdots + 1 / k_{r-2}}
  \Bigl(
    \int_{1/H}^{1/2} \vert \E_{k_{r-1}}(\alpha) \vert^2
      \, \frac{\dx \alpha}\alpha
    \int_{1/H}^{1/2} \vert \E_{k_r}    (\alpha) \vert^2
      \, \frac{\dx \alpha}\alpha
  \Bigr)^{1/2} \\
\label{bound-B-rh-minor}
  &\ll_{\K}
  N^{1 / k_1 + \cdots + 1 / k_{r-2}}
  \bigl(
    N^{1 / k_{r-1} + 1 / k_r} L^6
  \bigr)^{1/2}
  \ll_{\K}
  N^{\rho - 1 / (2 k_{r-1}) - 1 / (2 k_r)} L^3.
\end{align}
The other summands in $I_3$ can be treated in the same way, by
\eqref{E-bound} again.

Summing up from \eqref{final-mt-rh}, \eqref{bound-A-rh},
\eqref{bound-B-rh}, \eqref{bound-A-rh-minor},
\eqref{bound-B-rh-minor} and recalling that
$2 \le k_1 \le \dots \le k_r$, we proved that
\[
  \sum_{n = N + 1}^{N + H}
    \e^{-n / N} R(n; \K)
  =
  \frac{G(\K)}{\e \Gamma(\rho)} H N^{\rho - 1}
  +
  \Odip{\K}{\Phi_{\K}(N, H)},
\]
where, dropping terms that are smaller in view of the constraint
in~\eqref{first-bound-H}, we set
\begin{equation}
\label{def-Phi}
  \Phi_{\K}(N, H)
  =
  H^2 N^{\rho - 2}
  +
  H^{1 / 2} N^{\rho - 1/2 - 1 / (2 k_r)} L^3.
\end{equation}
Since we want an asymptotic formula, we need to impose the restriction
\begin{equation}
\label{bound-H-RH}
  H = \infty \bigl( N^{1 - 1 / k_r} L^6 \bigr),
\end{equation}
which supersedes \eqref{first-bound-H}.

We remark that when $k_1 = 2$ we can use Lemma~2
of~\cite{LanguascoZ2016c} instead of Lemma~\ref{Tolev-Lemma} in the
partial integration leading to \eqref{Tolev-appl}, and we can replace
the right-hand side by $N^{1 / 2} L^2 + H L^2$.
This means, in particular, that, in this case, we may replace $L^3$
in the far right of \eqref{def-Phi} by $L^{5 / 2}$.

Next, we remove the exponential weight, arguing essentially as in
\S\ref{final-Th1}.
This completes the proof of Theorem~\ref{Th-RH}.

\providecommand{\MR}{\relax\ifhmode\unskip\space\fi MR }
\providecommand{\MRhref}[2]{%
  \href{http://www.ams.org/mathscinet-getitem?mr=#1}{#2}
}
\providecommand{\href}[2]{#2}

\bigskip

\begin{tabular}{l}
Marco Cantarini \\
Dipartimento di Matematica e Informatica \\
Universit\`a di Perugia \\
Via Vanvitelli, 1 \\
06123, Perugia, Italia \\
email (MC): \texttt{marco.cantarini@unipg.it} \\
Alessandro Gambini, Alessandro Zaccagnini \\
Dipartimento di Scienze, Matematiche, Fisiche e Informatiche \\
Universit\`a di Parma \\
Parco Area delle Scienze 53/a \\
43124 Parma, Italia \\
email (AG): \texttt{a.gambini@unibo.it} \\
email (AZ): \texttt{alessandro.zaccagnini@unipr.it}
\end{tabular}

\end{document}